УДК 511.3

S.V. Matnyak

The asymptotic formula in the Goldbach binary problem

Khmelnytsky, Ukraine

Summary: The binary problem of Goldbach is solved by the method of the trigonometrical sums. The asymptotic formula of distribution of the even numbers formed by the sum of two simple uneven numbers is found for each even number from the set of natural numbers $N$.

Key words: binary problem of Goldbach, asymptotic formula, simple uneven numbers, even numbers, method of the trigonometrical sums.

**Introduction.** In this paper well-known unsolved Goldbach problem is considered. This problem represents itself the hypothesis according to which some even number, no less than six, is represented by the sum of two prime uneven numbers, which is the Goldbach binary problem .Goldbach had come to this problem experimentally. Nowadays computer gives you to experimentally test this hypothesis. But until the present moment no example which would contradict the Goldbach hypothesis has been found. A great progress in its proving was achieved by I.M. Vynogradov and L.G.Shnirelman. In 1973 Vynogradov managed to show that some sufficiently large uneven number is a sum of three prime uneven numbers. In this paper the asymptotic formula of distribution of sums of two prime uneven numbers is found, that is the founding of asymptotic formula faund gives the solution to the binary problem of Goldbach.

**Problem statement.** To find an asymptotic formula of distribution of even numberslarger than 6, which consist of addition of two prime uneven numbers.

**Solution**. In the proof presented here the method of I.M. Vynogradov is used, which was ingested to prove the ternary Goldbach problem.

**Definitions and theorems to reminder.**

**Definition 1.** The even number is a positive integer, if there is a number 2 among its prime multipliers.

**Definition 2.** Uneven number is a positive integer, if there no number 2 among its prime multipliers.

**Definition 3.** [1, p.119] The integer number p is called prime, if it differ from 0 and $\pm 1$ and has a divisor only to $\pm 1$ and $\pm p$. The integer number a, differ from 0 and $\pm 1$ and which has except $\pm 1$ and $\pm a$ other divisors, is called complex number.

**Theorem 1.** [1, p.120] Some integer positive number differs from 1, can be presented in the form of positive prime numbers product. This production accurate within multiplier order is unique.

**Theorem 2.** [1, p.120] The set of positive prime numbers is infinite.

**Theorem 3.** [1, p.394] (J. Hadamard and Vallée Poussin law). The ratio $\pi(x) : \frac{x}{\ln x}$ tends to 1 at unlimited growth $x$, i.e.

$$\lim_{x \to \infty} \frac{\pi(x) \cdot \ln x}{x} = 1,$$

where the function $\pi(x) = \frac{x}{\ln x}$ – defines the number of prime uneven numbers in the set of natural numbers $N$ when $x \to \infty$.

**Consequence.** (Euler's theorem) [1, p.408]. Each positive prime uneven number $a$ can be presented in the following form

$$a = 2 \cdot q + 1,$$

where $q$ – is an even number among the set of natural numbers $N$.

**Theorem 4.** The sum of arbitrary two prime and uneven numbers is an even number

$$N = p_1 + p_2,$$

Where $N$ is an even number, and $p_1, p_2$ are prime uneven numbers.

**Proof.** In terms of Euler's theorem consequence each uneven prime number $p_1$ and $p_2$ can be represented in following form

$$p_1 = b_1 \cdot q_1 + r_1 = 2q_1 + 1, \quad p_2 = b_2 \cdot q_2 + r_2 = 2 \cdot q_2 + 1$$

Where $q_1$ и $q_2$ – a whole, positive number from the set of natural numbers N.

Then

$$N = p_1 + p_2 = b_1 \cdot q_1 + r_1 + b_2 \cdot q_2 + r_2 = 2 \cdot q_1 + 1 + 2 \cdot q_2 + 1 = 4 \cdot (q_1 + q_2) + 2 = 2 \cdot (2 \cdot (q_1 + q_2) + 1).$$

Let us factorize the number $N$ is to prime numbers.

$$N = 2 \cdot m = 2 \cdot (2 \cdot (q_1 + q_2) + 1) = 2 \cdot p_1^{\alpha_1} \cdot p_2^{\alpha_2} \cdot \ldots \cdot p_s^{\alpha_s},$$

where numbers $p_1, p_2, \ldots, p_s$ are different positive prime numbers, and $\alpha_i > 0$
For $i = 1, 2, 3, \ldots, s$. Therefore, the number $N$ is even number, because between his simple multipliers is the number 2.

**Lemma 1.** (Paige). If for the given $\varepsilon_0$ quite high $c_1$ and $c$ are set. Then the number $\pi(N, q, l)$ of prime numbers, does not exceed $N$, which is in an arithmetic progression $q \cdot x + l$, $0 < q < r^c, (q, l) = 1, 0 < l < q$, $r = \ln N$, expressed by the formula

$$\pi(N, q, l) = \frac{1}{q} \int_2^N \frac{dx}{\ln x} + H, q_1 = \varphi(q),$$

where for all $q$, except on the series of different numbers, multiples of some $q = q_0$, the condition is met

$$q \geq r^{2-\varepsilon_0}$$

we have the inequality

$$H \ll \frac{N \cdot r^{-c}}{q_1 \cdot r}.$$

**Proof.** This Lemma is a consequence of the well known theorem of English mathematician Page.

**Lemma 2.** If $N \geq 2$, $r = \ln N$, $z$ – actual

$$I(z) = \int_2^N \frac{e^{2\pi i z x}}{r} dz, \quad J(z) = \int_2^N \frac{e^{2\pi i z x}}{\ln x} dx,$$

Then we have

$$I(z) \ll Z, \quad J(z) \ll Z,$$

$$Z = \begin{cases} N \cdot r^{-1}, & then |Z| \leq N^{-1} \\ |z|^{-1} \cdot r^{-1}, & then N^{-1} < |z| \leq N^{-0,5} \end{cases}.$$

**Proof**. This Lemma is proved by Vynogradov [2,p.95].

**Lemma 3.** If $\tau = N \cdot r^{-c}$, where $c \geq 2$, and if

$$R = \int_{-\tau^{-1}}^{\tau^{-1}} (J(z))^2 \cdot e^{-2\pi i z N} dz,$$

Where $J(z)$ has the meaning specified in Lemma 2, then

$$R = \frac{N}{r^2} + O\left(\frac{N}{r^3}\right).$$

**Proof.** Let us apply the notation of Lemma 2. The Integral $R$ is compared with integral

$$I(z) = \int_{-0,5}^{0,5} (I(z))^2 \cdot e^{-2\pi i z} dz,$$

And we find

$$R - R_0 = \int_{-\tau^{-1}}^{\tau^{-1}} \left((J(z))^2 - (I(z))^2\right) \cdot e^{-2\pi i z N} dz + \left(\int_{-\tau^{-1}}^{\tau^{-1}} (I(z))^2 \cdot e^{-2\pi i z N} dz - R_0\right).$$

Here, the first component of the right part, according to Lemma 2 and inequality

$$|J(z) - I(z)| < \int_{2}^{N} \left(\frac{1}{\ln x} - \frac{1}{r}\right) \cdot dx << \frac{N}{r^2},$$

So we have

$$<< \int_{-\tau^{-1}}^{\tau^{-1}} Z \cdot \frac{N}{r^2} dx << \int_{0}^{N^{-1}} \frac{N^2}{r^3} dz + \int_{N^{-1}}^{\tau^{-1}} \frac{N}{r^3} dz << \frac{N}{r^3}.$$

Another component of the right part, according to Lemma 2 will be

$$<< \int_{\tau^{-1}}^{0,5} Z^2 dz << \int_{\tau^{-1}}^{0,5} \frac{dz}{z^2 \cdot r^2} << \frac{N}{r^3}.$$

Therefore, $R - R_0 << \frac{N}{r^3}$. Let us assume,

$$R' = \int_{-0,5}^{0,5} (S(z))^2 \cdot e^{-2\pi i z N} dz, \quad S(z) = \sum_{x=2}^{N} \frac{e^{2\pi i z x}}{r},$$

In accordance to Lemma 2, we'll have $I(z) - S(z) << \frac{1}{r}$, $Z >> \frac{1}{r}$. Then

$$R_0 - R' << \int_0^{0,5} Z \cdot \frac{1}{r} dz = \int_0^{N-1} \frac{N^2}{r^2} dz + \int_{N-1}^{0,5} \frac{1}{r^2 \cdot z} dz \approx \frac{N}{r^2}.$$

Therefore, $R - R' << \frac{N}{r^2}$. The expression $r^2 \cdot R'$ defines the number of representations of number $N$ in the form of

$$N = x_1 + x_2$$

with integer $x_1$, $x_2$, which are greater than the number 2. The number of combinations of two numbers with $N-2$ will be equal to

$$\overset{2}{\underset{N-2}{C}} = \frac{(N-2) \cdot (N-3)}{2} = \frac{N^2 - 5N + 6}{2}.$$

Therefore, each number of a series of positive integer $N$ is possible to represent as

$$r^2 \cdot R' = \frac{\overset{2}{\underset{N-2}{C}}}{N} = \frac{N^2 - 5N + 6}{N} = \frac{N}{2} + O(1)$$

times, and each even number of series of positive integer $N$ is possible to represent in the form of

$$r^2 \cdot R' = \frac{2 \cdot \overset{2}{\underset{N-2}{C}}}{N} = \frac{2 \cdot (N^2 - 5N + 6)}{2N} = N + O(1)$$

times. From where we verify the validity of Lemma.

**Theorem 5.** Number $I(N)$ of representations of pair positive number $N$ as the sum

$$N = p_1 + p_2$$

simple unpaired numbers $p_1$ и $p_2$ and is expressed by the formula

$$I(N) = \frac{N}{2 \cdot r^2} \cdot S(N) + O\left(\frac{N}{r^{3-\varepsilon}}\right),$$

Where
$$S(N) = 3 \cdot \Pi_{(p,N)=1}\left(1 + \frac{1}{p^2}\right) \cdot \Pi''_{p/N}\left(1 - \frac{2}{p^2+1}\right)$$

More over $\Pi''_{p/N}$ applies to all of the prime numbers that divide $N$, a $\Pi_{(p,N)=1}$ -on prime divisors, which do not divide the number $N$ и $p < N$.

And we have $S(N) > 1,0$.

**Consequence (the binary Gold Bach problem).** There is $c_0$ with the condition that a certain pair number $N$ not less than $c_0$, is a sum of two prime's unpaired numbers.

**Proof.** Let us assume $\tau = N \cdot r^{-7}$, we have

$$I(N) = \int_{-\tau^{-1}}^{-\tau^{-1}+1} S_\alpha^2 \cdot e^{-2\pi i \alpha N} d\alpha, \; S_\alpha = \sum_{2 < p \leq N} e^{2\pi i \alpha p},$$

the interval of integration of the integral $I(N)$ can be divided into intervals of the first and second classes. Intervals of the first class we'll call intervals, which include all the values of alpha bud.

$$\alpha = \frac{a}{q} + z, (a,q) = 1, -\tau^{-1} < z < \tau^{-1}, 0 < q < r^2.$$

Obviously, ($c_0$ is big enough), intervals of first class do not intersect. Intervals of the second class are called the intervals that remained after the separation of the intervals of the first class. The arbitrary alpha, which belongs to the second class can be represented in the form

$$\alpha = \frac{a}{q} + z, (a,q) = 1, -\frac{1}{q\tau^{-1}} \leq z \leq \frac{1}{q\tau^{-1}}, r^2 < q \leq \tau.$$

Accordingly to the specified distribution of intervals the integration of the integral $I(N)$ is split into two components.

$$I(N) = I_1(N) + I_2(N).$$

1. **Assessment $I_2(N)$.** Accordingly to the theorem of the 1 chapter 4 by

$q \geq r^7$ we have

$$S_\alpha \ll N \cdot r \cdot \ln r \cdot \left( \sqrt{\frac{1}{q} + \frac{qr}{N}} + r \cdot e^{-0.5\sqrt{r}} \right) \ll \frac{N}{r^{2.5-\varepsilon_1}}.$$

And respectively to theorem 2 of the chapter 4 by $r^2 < q \leq r^7$, valid is the ratio is valid

$$S_\alpha \ll N \cdot r^{-2.5+\varepsilon_1}, \quad A = \frac{N}{[r]^{12}}.$$

Therefore

$$I_2(N) \ll N \cdot r^{-2.5+\varepsilon_1} \cdot \int_0^1 \left( \sum_{2<p'\leq N} \sum_{2<p<N} e^{2\pi i \alpha (p'-p)} \right)^{\frac{1}{2}} d\alpha \leq \frac{N}{r^{3-\varepsilon_1}}.$$

1. **Lemma 4.** Integral $\int_0^1 \left( \sum_{2<p'<N} \sum_{2<p<N} e^{2\pi i \alpha (p'-p)} \right)^{\frac{1}{2}} d\alpha \leq \frac{1}{\sqrt{\ln N}}$ under conditions, that $\alpha > \frac{\ln N}{2 \ln p}$ and $a > \frac{q}{2}$.

**Proof.** Let us write the expression $\int_0^1 \left( \sum_{2<p'<N} \sum_{2<p<N} e^{2\pi i \alpha (p'-p)} \right)^{\frac{1}{2}} d\alpha$ in the form of

$\int_0^1 \left| \sum_{2<p<N} e^{2\pi i \alpha \cdot \ln p} \right|^{\frac{1}{2}} d\alpha$, where $p' = p + \ln p$. Let us estimate the amount of

$\left| \sum_{2<p<N} e^{2\pi i \alpha \cdot \ln p} \right| \leq \left| \sum_{h=1}^{\left[\frac{N}{\ln N}\right]} e^{2\pi i h \eta} \right| = \left| \frac{e^{2\pi i \left( \left[\frac{N}{\ln N}\right]+1 \right) \eta} - e^{2\pi i \eta}}{e^{2\pi i \eta} - 1} \right| \leq \frac{2}{|e^{2\pi i \eta} - 1|} = \frac{1}{|\sin \pi \eta|} \leq \frac{1}{2\eta} \leq \frac{1}{\ln N}$, where

$\eta = \alpha \cdot \ln p$ and $\alpha > \frac{\ln N}{2 \ln p}$. From here we find, that $2\alpha > 1$ or $a > \frac{q}{2}$. From theorem 4 (Bertrand's postulate) [3, p. 97] we have, that $a$ should be in the interval $\left( \frac{N}{2} < a \leq N \right)$, that meets the condition $a > \frac{q}{2}$ as far as $1 < a < q < N$. Therefore

$$\int_0^1 \left( \sum_{2<p'<N} \sum_{2<p<N} e^{2\pi i \alpha (p'-p)} \right)^{\frac{1}{2}} d\alpha \leq \frac{1}{\sqrt{\ln N}}.$$

Lemma 4 is proved.

2. **The intervals of the first class, which correspond to the value of q which is not different.** If $I_{a,q}$ is the part of the integral that corresponds to the interval $I_1(N)$ of the first class, which includes quotient $\dfrac{a}{q}$ with value $q$, and not differ. Take the arbitrary $\alpha = \dfrac{a}{q} + z$, of this interval, the sum $S_\alpha$ can be divided into $\dfrac{1}{[r]^{12}}$ additive aspect

$$S_{\alpha,N_1} = \sum_{N_1-A<p\leq N_1} e^{2\pi i\left(\frac{a}{q}+z\right)p}, \quad A = \dfrac{N}{[r]^{12}},$$

For additive sums $S_{\alpha,N_1}$ we have $|z\cdot N_1 - z\cdot p| \leq z\cdot A$. Number of this summands $<< A\cdot r^{-1}$. Therefore

$$S_{\alpha,N_1} - e^{2\pi i z N_1}\sum_{N_1-A<p\leq N_1} e^{2\pi i\left(\frac{a}{q}+z\right)p} << z\cdot A^2\cdot r^{-1} << \dfrac{A}{r^6}.$$

But (Lemma 1)[2, p. 65] c=17) when given $l$ with the condition $0\leq l<q$, $(l,q)=1$, the number of primes of the form $q\cdot x + l$, which lie in the interval $N_1-A<p\leq N_1$ is expressed by the formula

$$\dfrac{1}{q_1}\int_{N_1-A}^{N_1}\dfrac{dx}{\ln x} + O\!\left(\dfrac{N\cdot r^{-17}}{q_1 r}\right)$$

Therefore:

$$S_{\alpha,N_1} = \sum_l e^{2\pi i\frac{a}{q}l}\dfrac{1}{q_1}\int_{N_1-A}^{N_1}\dfrac{e^{2\pi i N_1}dx}{\ln x} + O\!\left(\dfrac{A}{r^6}\right)$$

We shall find symbol (Lemma 2)

$$\sum_l e^{2\pi i\frac{a}{q}l} = \mu(q), \ |z\cdot N_1 - z\cdot p|\leq z\cdot A, \ \int_{N_1-A}^{N_1}\dfrac{zAdx}{\ln x} << \dfrac{A}{r^6},$$

$$S_{\alpha,N_1} = \dfrac{\mu(q)}{q_1}\int_{N_1-A}^{N_1}\dfrac{e^{2\pi i z x}dx}{\ln x} + O\!\left(\dfrac{A}{r^6}\right),$$

$$S_\alpha = \frac{\mu(q)}{q_1} J(z) + O\left(\frac{A}{r^6}\right),$$

$$\frac{J(z)}{q_1} \ll \frac{z}{q_1}, \quad S_\alpha^2 - \frac{\mu(q)}{q_1}(J(z))^2 \ll \frac{z}{q_1} Nr^{-6} + N^2 \cdot r^{-12}.$$

We denote $\tau_1 = \dfrac{N}{r^2}$ and because the range

$$\left(-\tau_1^{-1}, \tau_1^{-1}\right) \text{ enters into } the\ range\left(-\tau^1, \tau^1\right), \text{ so}$$

$-\tau_1^{-1} < z < \tau_1^{-1}$ and then

$$I_{\alpha,q} - \int_{-\tau_1^{-1}}^{\tau_1^{-1}} \frac{\mu(q)}{q_1^2}(J(z))^2 \cdot e^{2\pi i\left(\frac{a}{q}+z\right)N} dz \ll \int_0^{\tau_1^{-1}} \left(\frac{Z}{q_1} Nr^{-6} + \frac{N^2}{q_1} r^{-8}\right) dz \ll$$

$$\ll \frac{z}{q_1} \cdot N \cdot r^{-6} \frac{r^2}{N} + \frac{N^2}{q_1} \cdot r^{-8} \cdot \frac{r^2}{N} \ll \frac{N}{r^5 \cdot q_1^1},$$

Compel $a$ to runt his deduction system on module $q$, $therefore\ we\ obtain$

$$\sum_a I_{a,q} = G(q) \cdot R + O\left(N \cdot r^{-3} \cdot q_1^{-1}\right),$$

$$G(q) = \frac{\mu(q)}{q_1^2} \sum_a e^{2\pi i \frac{a}{q} N}, \tag{1}$$

According Lemma 3 and inequality $G(q) \ll q_1^1$, we shall write down

$$\sum_a I_{a,q} = \frac{N}{r^2} G(q) + O\left(\frac{N}{r^3 \cdot q_1}\right),$$

**3. Basic intervals, which correspond to the different values of q.** Let's assume that $q$ belongs to the different quantity of numbers $q$, So $q = q_0 \cdot k$ where $k$ is the number with condition $0 < k \leq r^{1+\varepsilon_0}$ (since $q \leq r^2$ and $q_0 \geq r^{1-\varepsilon_0}$).

dale us assume that

$$\alpha = \frac{a}{q} + z, -\tau \leq z \leq \tau,$$

$$\alpha = \frac{a}{q} + z, \quad -\tau \leq z \leq \tau, \quad \delta = |z|N,$$

$$S(q,\delta) = N \cdot r^{-1+\varepsilon} \cdot q^{-0,5}, \qquad \text{when} \delta \leq 1,$$

$$S(q,\delta) = N \cdot r^{-1+\varepsilon} \cdot q^{-0,5} \cdot \delta^{0,5}, \qquad \text{when} \delta \geq 1,$$

According to theorem 2 in Chapter 4[2], we will have $S_\alpha \ll S(q,\delta)$, So we obtain

$$I_{a,q} \ll \int_0^{\tau^{-1}} (S(q,\delta)^2) dz \ll \frac{1}{N}\int_0^1 (S(q,\delta)^2) d\delta \ll \frac{N}{r^{2-2\varepsilon}.kq_0}.$$

Compel $a$ to run this deduction system on module $q$, therefore we obtain

$$I_{a,q} \ll \frac{N}{r^{2-2\varepsilon}q_0}.$$

or, since the order of the first member is lower than of the following second order, so

$$\sum_a I_{a,q} = \frac{N}{r^2} G(q) + O\left(\frac{N}{r^{2-\varepsilon}r}\right).$$

**4.** The previous formulas $I(N)$. According to the computation (1,2,3) we find

$$I_1(N) - \sum_{q \leq r^2} \frac{N}{r^2}.G(q) \ll \sum_{q \leq r^2} \frac{N}{r^3 q_1} + \sum_{k \leq r^{1+\varepsilon}} \frac{N\cdot}{q_0^1 \cdot k^1 r^{3-2\varepsilon}} \ll \frac{N \cdot r^{\varepsilon_1}}{r^3},$$

$$\sum_{q>r^2} \frac{N}{r^2}.G(q) \ll \sum_{q>r^2} \frac{N}{r^2}.q_1^{-2} \ll \frac{N}{r^3},$$

$$I(N) = \frac{N}{r^2} \cdot S(N) + O\left(\frac{N}{r^{3-\varepsilon}}\right),$$

$$S(N) = \sum_{q=1}^\infty G(q). \qquad (2)$$

**5. Inverse transform and researching** $S(N)$. It is clear that, $G(n)$ differs from zero only in case, when cannonical expansion of the number of $q$ is given by $q = p_1 \cdot p_2 \cdot \ldots \cdot p_k$ (and in case, when $q=1$). Thus the identity will be true also if

$$G(p_1) \cdot \ldots \cdot G(p_k) = G(p_1 \cdot \ldots \cdot p_k)$$

Indeed, in the case $k=2$, the validity of this equation comes from the equality

$$G(p_1) \cdot G(p_2) = \frac{1}{\varphi(p_1 p_2)} \cdot \sum_{0 < a_1 < p_1} \sum_{0 < a_1 < p_2} e^{2\pi i \frac{(a_1 \cdot p_2 + a_2 \cdot p_1) \cdot N}{p_1 \cdot p_2}},$$

Where $a_1 \cdot p_2 + a_2 \cdot p_1$ runs this deduction system on module $p_1 p_2$. The generalization of this identity on the case $k > 2$ is trivial.

When $x \geq 2$ we have

$$\Pi_{p \leq x}(1 + G(p)) = \sum_{q \leq x} G(q) + \sum'_{q > x} G(q),$$

Where $\Sigma'$ extends to the values $q$, which are divided into primes, greater than $x$. As $S(N)$ converges absolutely at boundless increase $x$, the first components the right part aspires to a limit, the second to a limit zero. Therefore, forcing $p$ to run all simple numbers, we will have

$$S(N) = \Pi_3(1 + G(p)).$$

Let's consider the expression (1) $N$ is divided into on $q$, to $\sum_a e^{2\pi i \frac{a}{q} N} = q$. Let's substitute this expression in (1) and we will receive

$$G(q) = \frac{\mu(q) \cdot q}{q^2 \cdot \left(1 - \frac{1}{p_1}\right)^2 \cdot \ldots \cdot \left(1 - \frac{1}{p_k}\right)^2} = \frac{\mu(q)}{q \cdot \left(1 - \frac{1}{p_1}\right)^2 \cdot \ldots \cdot \left(1 - \frac{1}{p_k}\right)^2},$$

Where there is Euler's function which equals [4, item 30]

$$\varphi(q) = q \cdot \left(1 - \frac{1}{p_1}\right) \cdot \left(1 - \frac{1}{p_2}\right) \cdot \ldots \cdot \left(1 - \frac{1}{p_k}\right).$$

Let's substitute value in expression (2) $\left(S(N) = \sum_{p/N} G(q)\right)$ and obtain:

$$S(N) = \sum_{q/N} G(q) = \sum_{q/N} \frac{\mu(q)}{q \cdot \left(1 - \frac{1}{p_1}\right)^2 \cdot \ldots \cdot \left(1 - \frac{1}{p_k}\right)^2} = \frac{p_1^2}{(p_1 - 1)^2} \cdot \ldots \cdot \frac{p_k^2}{(p_k - 1)^2} \cdot \sum_{q/N} \frac{\mu(q)}{q} =$$

$$= \frac{p_1^2}{(p_1 - 1)^2} \cdot \ldots \cdot \frac{p_k^2}{(p_k - 1)^2} \cdot \left(1 - \frac{1}{p_1}\right) \cdot \left(1 - \frac{1}{p_2}\right) \cdot \ldots \cdot \left(1 - \frac{1}{p_k}\right) = \frac{p_1}{(p_1 - 1)} \cdot \frac{p_2}{(p_2 - 1)} \cdot \ldots \cdot \frac{p_k}{(p_k - 1)} > 1,$$

because [4,p.29]

$$\frac{p_1}{p_1 - 1} > 1, \frac{p_2}{p_2 - 1} > 1, \ldots, \frac{p_k}{p_k - 1} > 1 \text{ and than } \sum_{q/N} \frac{\mu(q)}{q} = \left(1 - \frac{1}{p_1}\right) \cdot \left(1 - \frac{1}{p_2}\right) \cdot \ldots \cdot \left(1 - \frac{1}{p_k}\right).$$

As $p_1 = 2$, we will have

$$S(N) > 2,0.$$

**and** $S(N) > 2,0 > \prod_{p \backslash N}\left(1 - \frac{1}{p^2}\right) = \frac{6}{\pi^2} > 0,6.$

**Then it is possible to write down that**

$$G(p) = -\frac{1}{p^2}, \text{ when } N \text{ divided on } p.$$

Therefore:

$$S(N) = 3 \cdot \Pi'\left(1 + \frac{1}{p^2}\right) \cdot \Pi'' \frac{\left(1 - \frac{1}{p^2}\right)}{\left(1 + \frac{1}{p^2}\right)} = 3 \cdot \Pi'\left(1 + \frac{1}{p^2}\right) \cdot \Pi''\left(1 - \frac{2}{p^2 + 1}\right),$$

where $\Pi'$ – extends on values $p$, which do not divide by $N$, and $\Pi''$ – on values which divide $N$ ...

Here the first multiplier $\Pi'$ of the right part is more than 1, and the second multiplier $\Pi''$ equals

$$\prod_{p}\left(1 - \frac{1}{p^2}\right) = \frac{6}{\pi^2} > 0,6.$$

Therefore we find

$$S(N) > 1,0$$

**The theorem is proved.**

On August, 06th, 2013
Remade
April, 22nd, 2014